\documentclass[12pt]{article}

\usepackage[reqno]{amsmath}
\usepackage{amssymb, amsthm}
\usepackage{enumerate,arydshln} 
\newtheoremstyle{plainsl}%
	{\topsep}
	{\topsep}
	{\slshape} 
	{}
	{\normalfont\bfseries}
	{.}
	{ }
	{}

\swapnumbers

\theoremstyle{plainsl}
\newtheorem{theorem}{Theorem}[section]
\newtheorem{lemma}[theorem]{Lemma}
\newtheorem{corollary}[theorem]{Corollary}

\newtheorem{conj}[theorem]{Conjecture}

\renewcommand\proof{\noindent\textsl{Proof. }}
\newcommand\sqr[2]{{\vbox{\hrule height.#2pt
    \hbox{\vrule width.#2pt height#1pt \kern#1pt
        \vrule width.#2pt}\hrule height.#2pt}}}

\renewcommand\qed{%
	\ifmmode\eqno\sqr53
	\else\nolinebreak\ \hfill\sqr53\medbreak\fi}

\numberwithin{equation}{section}

\newcommand{\remove}[1]{}

\newcommand\sto{{\scriptsize $\to$}}

\newcommand\cA{{\mathcal A}}
\newcommand\seq[3]{#1_{#2},\ldots,#1_{#3}}
\newcommand\la{\lambda}
\newcommand\one{{\bf1}}

\title{A new proof of the Erd\H{o}s-Ko-Rado theorem for intersecting families of permutations}
\author{Chris Godsil \footnote{Research supported by NSERC.}\\
\small  Department of Combinatorics and Optimization \\[-0.8ex]
\small University of Waterloo,  Waterloo, Ontario, Canada\\[-0.8ex]
\small \texttt{cgodsil@math.uwaterloo.ca}\\
Karen Meagher \small{*}\\ [-0.8ex]
\small  Department of Mathematics and Statistics \\[-0.8ex]
\small University of Regina,  Regina, Saskatchewan, Canada\\[-0.8ex]
\small \texttt{kmeagher@math.uregina.ca}
}

\begin{document}
\maketitle

\abstract{Let $S(n)$ be the symmetric group on $n$ points. A subset
  $S$ of $S(n)$ is \textsl{intersecting} if for any pair of
  permutations $\pi, \sigma$ in $S$ there is a point $i \in
  \{1,\dots,n\}$ such that $\pi(i)=\sigma(i)$. Deza and
  Frankl~\cite{MR0439648} proved that if $S \subseteq S(n)$ is
  intersecting then $|S| \leq (n-1)!$.  Further, Cameron and
  Ku~\cite{MR2009400} show that the only sets that meet this bound
  are the cosets of a stabilizer of a point. In this paper we give a
  very different proof of this same result.}

\section{Introduction}

Cameron and Ku~\cite{MR2009400} proved a version of the
Erd\H{o}s-Ko-Rado theorem for permutations. In this paper we give an
alternate proof to this theorem which is substantially different from
the one given by Ku and Cameron.

The Erd\H{o}s-Ko-Rado theorem~\cite{MR25:3839} is a central result in
extremal combinatorics. There are many interesting proofs and
extensions of this theorem, for a summary see~\cite{MR86a:05004}. The
Erd\H{o}s-Ko-Rado theorem gives a bound on the size of a family of
intersecting $k$-subsets of a set and describes exactly which families
meet this bound.

\begin{theorem}(Erd\H{o}s, Ko and Rado~\cite{MR25:3839})
Let $k,n$ be positive integers with $n > 2k$.
Let $\cA$ be a family of $k$-subsets of $\{1,\dots,n\}$ such that any two sets 
from $\cA$ have non-trivial intersection, then $|\cA| \leq (n-1)!$.
Moreover, $|\cA| = (n-1)!$ if and only if $\cA$ is the collection of all $k$-subsets 
that contain a fixed $i \in \{1,\dots,n\}$.
\end{theorem}

The Erd\H{o}s-Ko-Rado theorem has been extended to objects other than
subsets of a set.  For example, Hsieh~\cite{MR0382015} and Frankl and
Wilson~\cite{MR867648} give a version for intersecting subspaces of a
vector space over a finite field, Berge~\cite{MR0389636} proves it for intersecting
integer sequences, Rands~\cite{MR84i:05024} extends it to intersecting
blocks in a design and Meagher and Moura~\cite{MR2156694} prove a
version for partitions.

The extension we give here is to intersecting permutations.  Let
$S(n)$ be the symmetric group on $\{1,\dots,n\}$. Permutations $\pi,
\sigma \in S(n)$ are said to be \textsl{intersecting} if $\pi(i) =
\sigma(i)$ for some $i \in \{1,\dots,n\}$.  Similar to the case for
subsets of a set, there are obvious candidates for maximum
intersecting systems of permutations, these are the sets
\begin{eqnarray}\label{eq:maxindy}
S_{i,j} = \{ \pi \in S(n) : \pi(i)=j \}, \quad i,j \in \{1,\dots,n\}.
\end{eqnarray}
These sets are the cosets of a stabiliser of a point.

\begin{theorem}(Cameron and Ku~\cite{MR2009400})\label{thm:main}
Let $n\geq 2$. If $S \subseteq S(n)$ is an intersecting family of permutations then:
\begin{enumerate}[(a)]
\item $|S| \leq (n-1)!$.
\item if $|S|=(n-1)!$ then $S$ is a coset of a stabiliser of a point.
\end{enumerate}
\end{theorem}

The proof given by Cameron and Ku uses an operation called
\textsl{fixing} which is similar to the shifting operation used in the
original proof of Erd\H{o}s-Ko-Rado. They show that a maximum
intersecting family of permutations is closed under this fixing
operation.  Assuming that the family contains the identity
permutation, and thus each permutation in the family has a fixed
point, they next consider the set system formed by the sets of fixed points
for each permutation in the family.  Cameron and Ku prove that if the
family of permutations is closed under the fixing operation, then this
set system is an intersecting set system. Finally, they prove the
result by showing that if a family of intersecting permutations has
size $(n-1)!$, then the sets in the intersecting set system must all
intersect in the same point.

Our proof uses a graph called the \textsl{permutation graph} which
appears in the paper by Cameron and Ku. This graph is a union of
graphs in an association scheme, we use properties of this association
scheme together with information about the group representation of the
symmetric group to get the result.

This approach has been used to prove the standard Erd\H{o}s-Ko-Rado theorem for
sets~\cite[Section 5.4]{mikethesis} and also to prove versions of the
Erd\H{o}s-Ko-Rado theorem for other objects such as the $3 \times 3$
uniform partitions and vector spaces over a finite
field~\cite{MR2260847}. It is interesting that this method also works
for permutations and hoped that this method can be generalized to
other objects.

The proof we give only applies for $n>6$, for smaller $n$ the result
can be verified using GAP~\cite{GAP4}.

\section{The Clique-Coclique Bound}

In this section we give a proof of the clique-coclique bound for the
union of graphs in an association scheme. Although this bound is not
new, it was originally proven by Delsarte~\cite{MR0384310}, and an
alternate proof for vertex-transitive graphs is given by Cameron and
Ku~\cite{MR2009400}, the proof given here is new.

Let $\cA=\{\seq A0d\}$ be an association scheme with $d$ classes on $v$ vertices
and let $v_i$ be the valency of the $i$-th graph.  
Denote the principal matrix idempotents of the association scheme
by $\seq E0d$ and let $m_i$ be the dimension of the eigenspace belonging to $E_i$.
We note that 
\[
E_0 =\frac1v J
\]
where $J$ is the all-ones matrix.

\begin{theorem}(Delsarte~\cite[Theorem 3.9]{MR0384310})
\label{thm:ccl}
Let $\cA$ be an association scheme on $v$ vertices and let $X$ be the
union of some of the graphs in the scheme.  If $C$ is a clique and $S$
is an independent set in $X$, then
\begin{equation}\label{ineq:cliquecoclique}
|C|\,|S| \le v.
\end{equation}
If equality holds and
$x$ and $y$ are the respective characteristic vectors of $C$ and $S$, then
\[
x^TE_jx\, y^TE_jy =0 \quad \textrm{for all } j>0.
\]
\end{theorem}

\proof
We have the following fundamental identity (see~\cite[Section 12.6]{MR1220704}):
\[
\sum_{i=0}^d\frac1{vv_i}x^TA_ix\, A_i = \sum_{j=0}^d \frac1{m_j}x^TE_jx\, E_j
\]
from which it follows that
\begin{equation}
\label{xaya}
\sum_{i=0}^d\frac1{vv_i}x^TA_ix\, y^TA_iy 
	= \sum_{j=0}^d \frac1{m_j}x^TE_jx\, y^TE_jy.
\end{equation}

Now suppose $C$ is a clique and $S$ is an independent set in $X$, and let $x$ and 
$y$ be their respective characteristic vectors.
The graph $X$ is a union of graphs in the scheme, if $A_i$ is one of the graphs in this union then $A_iy=0$ otherwise $A_ix=0$.
So for all $i>0$,
\[
x^TA_ix\,y^TA_iy=0,
\]
and hence the left side of Equation~\eqref{xaya} is
\begin{equation}
\label{xxyy}
\frac1{v}x^Tx\,y^Ty =\frac{|C|\,|S|}{v}.
\end{equation}
For all $j$, the matrix $E_j$ is positive semidefinite and therefore
\[
x^TE_jx\, y^TE_jy \ge0.
\]
Consequently the right side of Equation~\eqref{xaya} is bounded below by its first term:
\begin{equation}
\label{xex}
x^TE_0x\, y^TE_0y =\frac{1}{v^2}x^TJx\, y^TJy =\frac{|C|^2|S|^2}{v^2}.
\end{equation}
It follows from \eqref{xxyy} and \eqref{xex} that $|C|\,|S|\le v$, as required.
If equality holds the remaining condition follows immediately.\qed

We will prove a simple, but useful corollary of this result.

\begin{corollary}\label{cor:ccleql}  
Let $X$ be a union of graphs in an association scheme with the
property that the clique-coclique bound holds with equality. Assume
that $C$ is a maximum clique and $S$ is a maximum independent set in
$X$ with characteristic vectors $x$ and $y$ respectively.  If $E_j$
are the idempotents of the association scheme, then for $j>0$ at most
one of the vectors $E_jx$ and $E_jy$ is not zero.
\end{corollary}

\proof
If $j>0$, then
\[
x^TE_jx\, y^TE_jy =0.
\]
Since $E_j$ is positive semidefinite, $z^TE_jz=0$ if and only if $E_jz=0$.
\qed

\section{The Permutation Graph}

For a positive integer $n$ define the \textsl{permutation graph}
$P(n)$ to be the graph whose vertex set is the set of all permutations
of an $n$-set and vertices $\pi$ and $\sigma$ are adjacent if and only
if they are not intersecting, that is $\pi(i) \neq \sigma(i)$ for all
$i \in \{1, \dots ,n\}$.  The intersecting families of permutations
are exactly the independent sets in $P(n)$. We will show that the size
of the maximum independent set in $P(n)$ is $(n-1)!$ and the only sets
that meet this bound are the sets $S_{i,j}$ from
Equation~\ref{eq:maxindy}.

Let $d(n)$ be the number of derangements of an $n$-set (that is the
number permutation with no fixed points), then the graph $P(n)$ is
$d(n)$-regular. The number of derangements of a set of size $n$ is
defined by the following recursive formula
\begin{align}\label{eq:derangements}
d(n)=(n-1)\left( d(n-1)+d(n-2)\right)
\end{align}
with $d(1)=0$ and $d(2)=1$.

The permutation graph is a vertex-transitive graph, in fact, $P(n)$ is
a Cayley graph whose connection set is the set of all derangements.
Since this set is closed under conjugation, $P(n)$ is a
\textsl{normal} Cayley graph (for more on normal Cayley graphs
see~\cite[Section 5.2]{MR1468789}).

Further, the graph $P(n)$ is a union of graphs in the association
scheme known as the \textsl{conjugacy class scheme} on $S(n)$. The
conjugacy class scheme can be constructed for any group $G$ and is an
association scheme on the elements of $G$. Using the regular
representation each element of $G$ can be expressed as a $|G|
\times |G|$ permutation matrix. For any conjugacy class $C$ in $G$
define $A_C$ to be the sum of the permutation matrices for all the
elements in the conjugacy class. Then
\[
\mathcal{A} = \{A_C: C \textrm{ a conjugacy class in } G\}
\]
is the conjugacy class scheme on $G$ (for more on the conjugacy class
scheme see~\cite[page 54]{MR882540}).

If $\mathcal{A}$ is the conjugacy class scheme for the symmetric group
$S(n)$, then the adjacency matrix of $P(n)$ is the sum of $A_{C}$ over
all conjugacy classes $C$ of derangements. Since $P(n)$ is the sum of
graphs in an association scheme the clique-coclique bound
(Inequality~\ref{ineq:cliquecoclique}) holds. With this bound, it is
straightforward to get the first statement of Theorem~\ref{thm:main}.
This proof of the bound in Theorem~\ref{thm:main} is included
in~\cite[Theorem 5]{MR2009400} and it was also shown by Deza and
Frankl~\cite{MR0439648}.

\begin{theorem}\label{thm:clique}
The size of a maximum clique in $P(n)$ is n. 
\end{theorem}
\proof A clique in $P(n)$ can have no more than $n$ vertices. This is
clear since the image of 1 (or any other element in $\{1,\dots,n\}$)
must be distinct for each permutation in the clique.  Further, each
row of a Latin square of order $n$ is a permutation in $S_n$ and the
set of all rows in a Latin square of order $n$ is a clique of size $n$
in $P(n)$.  Since a Latin square of order $n$ exists for every $n$ the
theorem holds.  \qed

\begin{theorem}
The size of a maximum independent set in $P(n)$ is $(n-1)!$.
\end{theorem}
\proof
Since the graph $P(n)$ is a union of graphs in an association scheme the clique-coclique bound
holds for $P(n)$, that is
\[
\alpha(P(n)) \leq \frac{|V(P(n))|}{\omega(P(n))}.
\]
{}From Theorem~\ref{thm:clique}, $\omega(P(n)) =n$
and hence
\[
\alpha(P(n)) \leq (n-1)!.
\]
Finally, the sets $S_{i,j}$ from Equation~\ref{eq:maxindy} are
independent sets of size $(n-1)!$.  \qed

\section{Eigenvalues of $P(n)$}

In this section we will find two eigenvalues of the adjacency matrix
of $P(n)$. Eigenvalues of the adjacency matrix of $P(n)$ will simply
be refer to as the eigenvalues of $P(n)$.

\begin{lemma}
For all positive integers $n$ 
\[d(n) \quad \textrm{ and } \quad -\frac{d(n)}{n-1} \]
are eigenvalues for $P(n)$.
\end{lemma}
\proof 
Consider the independent set $S_{n,n}$ as defined in Equation~\ref{eq:maxindy}. The partition 
\[
\{S_{n,n}, V(P(n)) \setminus S_{n,n}\}
\] 
is the orbit partition of $S(1) \times S(n-1)$ acting on the vertices
of $P(n)$, hence it is an equitable partition.  The quotient graph of
$P(n)$ with respect to this partition is
\[
\left(    
\begin{array}{cc}
0 & d(n) \\ \frac{d(n)}{n-1} & d(n) - \frac{d(n)}{n-1}
\end{array}
\right).
\]
The eigenvalues of this quotient graph are $d(n)$ and 
$-\frac{d(n)}{n-1}$. 
Since the partition is equitable these are also eigenvalues for the
graph $P(n)$.
\qed

Since $P(n)$ is a $d(n)$-regular graph, $d(n)$ is the largest eigenvalue of $P(n)$. 
By Equation~\ref{eq:derangements}
\[
-\frac{d(n)}{n-1} = -(d(n-1)+d(n-2))
\] 
so this eigenvalue is also an integer.

The eigenvalues of a graph can be used to find bounds on the size of
the maximum independent sets. In particular, if $X$ is a $d$-regular
vertex-transitive graph with least eigenvalue $\tau$ then
\[
\alpha(X) \leq \frac{|V(X)|}{1-\frac{d}{\tau}}.
\]
This is known as the \textsl{ratio bound for independent sets}
(see~\cite[Lemma 9.6.2]{MR1829620} for a proof).  Ku~\cite{MR2302532}
has conjectured that the least eigenvalue of $P(n)$ is
$-\frac{d(n)}{n-1}$.  If this is true, then the ratio bound gives the
first part of Theorem~\ref{thm:main}.

The eigenvalues of a graph in a conjugacy class scheme and the
idempotents of the conjugacy class scheme can be determined by the
character table of the group. We will state these formulas for
the conjugacy class scheme on the symmetric group.

It is well-known that each irreducible character of $S(n)$ corresponds
to an integer partition of $n$. To denote that $\la$ is an integer
partition of $n$, we write $\la \vdash n$. If $\lambda \vdash n$, we
will represent the character of $S_n$ corresponding to $\lambda$ by
$\chi_\la$.  Each partition $\la$ of $n$ corresponds to a module, we
will call this the $\la$-module.  For more on the representation
theory of the symmetric group see~\cite[Chapter 4]{MR1153249}.

For each $\la \vdash n$ there is a principal idempotent in the scheme.
This idempotent is the $n! \times n!$ matrix whose entries are given by
\begin{align}\label{eq:proj}
(E_\la)_{\pi,\sigma} = \frac{\chi_\la(1)}{n!} \chi_\la(\pi^{-1}\sigma)
\end{align}
where $\pi,\sigma \in S(n)$.

For $C$ a conjugacy class in $S(n)$ the
eigenvalues of $A_C$ are 
\[
p_C^{\,\lambda} = \frac{|C|}{\chi_\lambda(1)} \chi_\lambda (c),  \quad c \in C
\]
where $\lambda$ ranges over all partitions of $n$ (for a proof of this
see~\cite[Chapter II, Section 2.7]{MR882540}).

It follows from this that the eigenvalues of $P(n)$ are
\[
\sum_{C}p_C^{\,\lambda}, \quad \lambda \vdash n
\]
where the sum is taken over all conjugacy classes of derangements.

For the partition $\la = [n]$ the value of $p_C^{\,[n]}$ is $|C|$ and thus
\[
 \sum_{C} p_C^{\,[n]} = \sum_{C} |C| = d(n)
\]
where the sum is taken over all conjugacy classes of derangements.

For any $x \in S(n)$ the value of $\chi_{[n-1,1]}(x)$ is one less than the number of fixed points in $x$, so for $C$ any conjugacy class of derangements $p_C^{\,[n-1,1]} = -\frac{|C|}{n-1}$.
Thus 
\[
\sum_{C} p_{C}^{\,[n-1,1]} = \sum_{C} -\frac{|C|}{n-1} = -\frac{d(n)}{n-1} 
\]
again, the sum is taken over all conjugacy classes of derangements.

\section{The $(n-1)$-module}

For a subset $S \subseteq S(n)$ let $v_S$ be the characteristic vector
of $S$ and if $S$ is one of the independent sets $S_{i,j}$ defined in
Equation~\ref{eq:maxindy}, then we will simply denote $v_S$ by
$v_{i,j}$. Throughout this section $\one$ will denote the all-ones
vector, the length of $\one$ will be clear from context.

We will first show that for any maximum independent set $S$ the vector
$v_S-\frac{1}{n}\one $ is in the module corresponding to the
representation $[n-1,1]$.  The next step will be to prove that the
vectors $v_{i,j}-\frac{1}{n}\one$ span the $[n-1,1]$-module.  Finally
we show that any characteristic vector of a maximum independent
set that is in this span must be one of $v_{i,j}$ for $i,j \in \{1,\dots,n\}$.

\begin{lemma}\label{lem:module}
Let $n$ be an integer with $n>6$.  Let $S$ be a maximum independent
set in $P(n)$ and $v_S$ be the characteristic vector of $S$.  Then the
vector $v_S-\frac{1}{n}\one$ is in the $[n-1,1]$-module.
\end{lemma}
\proof First, a simple calculation shows that $v_s-\frac{1}{n}\one$
is orthogonal to $\one$, so this vector is not in the $[n]$-module.

For $\la$ an integer partition of $n$ let $\chi_\la$ be the character of $S_n$
corresponding to $\la$. For $C$ a maximum clique in $P(n)$ define
\[
\chi_\lambda(C) = \sum_{x \in C} \chi_\lambda(x).
\]

If $\chi_\lambda(C) \neq 0$, then by Equation~\ref{eq:proj} $E_\la v_C
\neq 0$. By Corollary~\ref{cor:ccleql}, this implies at $E_\la v_S =
0$ which in turn implies that $E_\la (v_S -\frac{1}{n}\one)= 0$
for all partitions $\la \neq [n]$.  This
means that the vector $v_S -\frac{1}{n}\one$ is orthogonal to the
$\la$-module.  If this is true for every partition $\la \vdash n$
except $[n-1,1]$, then for every maximum independent set $S$ the vector
$v_S-\frac{1}{n}\one$ is in the $[n-1,1]$-module.  To prove this
theorem, we will show for every $\la \vdash n$ with $\la \neq [n-1,1]$
there is a maximum clique $C$ such that $\chi_\la (C) \neq 0$.

For $n>6$ there is a decomposition of the complete digraph on $n$
vertices into $n-1$ directed cycles~\cite{MR1986837}.  Each of these
directed cycles is a cycle of length $n$ in $S_n$. Moreover, no two
cycles in the decomposition share an edge so these cycles are adjacent
in $P(n)$.  Let $T$ be the $n$-clique whose elements are the
$n$-cycles in this decomposition together with the identity of $S(n)$.

Since every $x \in T$, except the identity, is an $n$-cycle for every $\la \vdash n$ the value of $\chi_\lambda(x)$ is the same. Thus
\begin{eqnarray*}
\chi_\lambda(T)&=& \sum_{x \in T} \chi_\lambda(x) \\
         &=& \chi_\lambda(1) + (n-1)\chi_\lambda(x) \quad \textrm{ $x$ an $n$-cycle.}
\end{eqnarray*}
Further, $\chi_\lambda(x) = \pm 1$ for every character $\chi_\lambda$
(for a proof of this see~\cite{MR1093239}). Since
$\chi_\lambda(1)$ is positive, if $\chi_\lambda(C)=0$, then
$\chi_\lambda(x) =-1$ and $\chi_\lambda(1) = n-1$.  For $n>6$ the only
partitions of $n$ with $\chi_\lambda(1) = n-1$ are $[n-1,1]$ and
$[2,1^{n-2}]$.

If $n$ is even, then for $x$ an $n$-cycle
$\chi_{[2,1^{n-2}]}(x) = 1$ so $\lambda$ must be $[n-1,1]$.

Finally, if $n$ is odd we need to prove that $\lambda$ is $[n-1,1]$.
To do this we construct a clique $T$ with $\chi_{[2,1^{n-2}]}(T) \neq
0$.  Consider an $n \times n$ Latin square with the first row
$(1,2,\dots ,n)$ and the second row $(2,1,n,3,4, \dots, n-1)$. Such a
Latin square exists since any Latin rectangle can be extended to a
Latin square~\cite{MR0013111}. The rows of this Latin square will the be permutations
in our clique. The first row corresponds to the identity permutation,
the second to an odd permutation.  The first row will contribute $n-1$
to the sum $\chi_{[2,1^{n-2}]}(T)$ and the second row will contribute
1.  Each of the last $n-2$ permutations will contribute no less than
$-1$ to the sum so the sum cannot be 0.
\qed

Next we give a basis for the $[n-1,1]$-module.

\begin{lemma}\label{lem:basis}
For any $i,j \in \{1, \dots ,n-1\}$ let $v_{i,j}$ denote the
characteristic vector of the independent set $S_{i,j} =\{\pi \in S(n) :
\pi(i)=j\}$.  The vectors $v_{i,j}-\frac{1}{n}\one$ form a
basis for the $[n-1,1]$-module.
\end{lemma}
\proof From Lemma~\ref{lem:module}, the vectors
$v_{i,j}-\frac{1}{n}\one$ are elements in the $[n-1,1]$-module. The
dimension of the $[n-1,1]$-module is $(n-1)^2$, so we only need to show
that these vectors are linearly independent.  Since $\one \not\in
\mathrm{span} \{v_{i,j}: i,j \in\{1,\dots,n-1\} \}$, it is enough to
show that the vectors $v_{i,j}$ are linearly independent.

Order the pairs in $\{1, \dots ,n-1\}$ so that pair $(i,j)$ occurs before $(k,\ell)$
if $i < k$ or if $i=k$ and $j < \ell$.  
Let $H$ be a $01$-matrix with size $n! \times (n-1)^2$ defined as
follows: the columns are indexed by the pairs from the
$(n-1)$-set in the above ordering and the rows are indexed by all the
permutations of an $n$-set. The $(\pi, (i,j))$-entry of $H$ is 1 if and
only if $\pi(i)=j$.

Let $I_{n}$ be the $n\times n$ identity matrix and $J_n$ the $n\times n$ all-ones matrix.
The adjacency matrix of the complete graph on $n$ vertices is $K_n = J_n -I_n$. 
It is not hard to see with the given ordering on the pairs that
\[
H^T H = (n-1)!I_{(n-1)^2} + (n-2)! (K_{n-1}\otimes K_{n-1}).
\]
Since 0 is not an eigenvalue of this matrix, it has rank $(n-1)^2$.
Finally, the rank of $H$ is equal to the rank of $H^T H$ and the result holds.
\qed

\section{Proof of Theorem~\ref{thm:main}}

Let $H$ be the $n! \times (n-1)^2$ matrix whose rows are the elements
of the symmetric group on $n$ points and columns are the ordered pairs
from $\{1,\dots,n-1\}$ with the $(\pi, (i,j))$ position of $H$ equal
to 1 if $\pi(i)=j$ and zero otherwise.

Denote the columns of $H$ by $h_{i,j}$. By Lemma~\ref{lem:basis} the vectors
\[
h_{i,j} -\frac{1}{n}\one, \quad i,j \leq n-1
\]
are a basis for the $[n-1,1]$-module.  By Lemma~\ref{lem:module}, for any independent set $S$, 
the vector $v_S - \frac{1}{n}\one$ is
in the $[n-1,1]$-module. In particular, it is in
\[
\mathrm{span}\left\{ h_{i,j} -\frac{1}{n}\one : i,j\in \{1,\dots,n-1\} \right\}.
\]
 
This implies that the characteristic vector of any maximal
independent set is in the span of column space of $H$ and $\one$.

Let $\sigma$ be the identity permutation on the $n$-set and let
$N(\sigma)$ denote the set of permutations adjacent to $\sigma$ in
$P(n)$ (these are the derangements). Consider three submatrices of
$H$:
\begin{enumerate}[(a)]
\item $N$ the submatrix whose rows are the permutations in $N(\sigma)$,
\item $M$ the submatrix of $N$ whose columns are all the pairs $(i,j)$
  with $i,j \in \{1,\dots,n-1\}$ and $i \neq j$,
\item $W$ the submatrix of $H$ whose columns are all the pairs $(i,i)$
  with $i\in \{ 1,\dots,n-1\}$.
\end{enumerate}

If the columns of $H$ are arranged so that the first $n-1$ columns
correspond to the pairs $(i,i)$ for $i =1,\dots,n-1$, and the rows are
arranged so the first row corresponds to the permutation $\sigma$ and
the next $d(n)$ rows correspond to the neighbours of $\sigma$, then
$H$ has the following block structure:
\begin{center}
\begin{tabular}{|c|c|} \hline
1 & 0 \\ \hline
0 & $M$ \\ \hline
$H_1$ & $H_2$  \\ \hline
\end{tabular} 
\end{center}
and the first $n-1$ columns form the matrix $W$.

\begin{lemma}
For all $n$ the rank of $M$ is $(n-1)(n-2)$.
\end{lemma}
\proof 
The matrix $K_{n-1} \otimes I_{n-2}$ has rank $(n-1)(n-1)$; we show it is a submatrix of $H$.
To find this submatrix, we reorder the rows and columns of $H$.

Order the pairs from $\{1,\dots,n-1\}$ so that the pair $(i,  i+j
\pmod{n-1})$ occurs before $(k, k+\ell\pmod{n-1})$ if $i < k$ or $i=k$
and $j < \ell$.  Order the columns of $M$ with this ordering.

Next we define an ordering on a subset of derangements.
Let $a \in \{1, \dots ,n-1\}$ and $b \in \{1, \dots ,n-2\}$. Define a
permutation of $\{1, \dots ,n\}$ for $i\in \{1, \dots ,n-1\}$ as follows:
\[
\pi_{a,b}(i) = \left\{  
    	\begin{array}{ll}
		n & \textrm{if } a = i; \\
		i+b & \textrm{if } a \neq i \textrm{ and } i+b < n; \\
		i+b+1 \pmod{n} & \textrm{if }  a \neq i \textrm{ and } i+b \geq n.
	\end{array}
		\right.
\]
Note that the value of $\pi_{a,b}(n)$ is forced.

Order these permutations so that $\pi_{a_1,b_1}$ occurs before
$\pi_{a_2,b_2}$ if $a_1 < a_2$ or $a_1 = a_2$ and $b_1 < b_2$.
Consider the submatrix of $M$ induced by the rows corresponding to the
permutations $\pi_{a,b}$ for $a \in \{1, \dots ,n-1\}$ and $b \in \{1,
\dots ,n-2\}$. This submatrix of $M$ is $K_{n-1} \otimes I_{n-2}$.
\qed

\begin{table}
\begin{center}
\begin{tabular}{cc cc:cc:cc}
$(a,b)$& $\pi_{a,b}$ & 1\sto 2 & 1\sto 3 & 2\sto 3 & 2\sto 1 & 3\sto 1 & 3\sto 2 \\ \cline{3-8}
(1,1) & (1,4,2,3)&\multicolumn{1}{|c} 0&0 & 1&0 & 1&\multicolumn{1}{c|} 0 \\
(1,2) & (1,4,3,2)&\multicolumn{1}{|c} 0&0 & 0&1 & 0& \multicolumn{1}{c|}1 \\
\cdashline{3-8}
(2,1) & (1,2,4,3)&\multicolumn{1}{|c} 1&0 & 0&0 & 1& \multicolumn{1}{c|}0 \\
(2,2) & (1,3,2,4)&\multicolumn{1}{|c} 0&1 & 0&0 & 0& \multicolumn{1}{c|}1 \\
\cdashline{3-8}
(3,1) & (1,2,3,4)&\multicolumn{1}{|c} 1&0 & 1&0 & 0& \multicolumn{1}{c|}0 \\
(3,2) & (1,3,4,2)&\multicolumn{1}{|c} 0&1 & 0&1 & 0& \multicolumn{1}{c|}0 \\
\cline{3-8}
\end{tabular}
\end{center}
\caption{The submatrix of $M$ for $n=4$.}
\end{table}

\begin{lemma}\label{lem:kernel}
If $y$ is in the kernel of $N$, then $Hy$ lies in the column space of $W$.
\end{lemma}
\proof Assume $y$ is in the kernel of $N$.  Let $y_M$
denote the vector of length $(n-1)(n-2)$ formed by taking the final
$(n-1)(n-2)$ entries of $y$. Then
\[
0 = Ny = [0| M ]y = My_M. 
\]
Since $M$ has rank $(n-1)(n-2)$, the last $(n-1)(n-2)$ entries of $y$ are all 0.
Thus $Hy$ is in the column space of $W$.
\qed

Let $[N|1]$ be the $d(n) \times ((n-1)^2+1)$ matrix with a column of
ones added to $N$, and $[M|1]$ the $d(n) \times ((n-1)(n-2)+1)$ matrix
with a column of ones added to $M$. As above, for a length $(n-1)^2+1$
vector $y$, the vector formed by the last $(n-1)(n-2)+1$ entries of
$y$ will be denoted by $y_{[M|1]}$.

\begin{lemma}\label{lem:kernelwith1}
If $y$ is in the kernel of $[N|1]$, then $y_{[M|1]}$ is a scalar multiple of
\[
(1,1,\dots,1,-(n-2)).
\]
\end{lemma}
\proof
As in the previous lemma,
\[
0 = [N|1]y = [0|M|1]y = [M|1]y_{[M|1]}.
\]
Since $M$ has full column rank, the dimension of the kernel of $[M|1]$
is at most 1.  Each row of $M$ has exactly $n-2$ entries equal to one
and all other entries zero, so the vector $(1,1,\dots,1,-(n-2))$ is in
the kernel of $[M|1]$ and is a basis for the kernel of $[M|1]$.  \qed

We now have all the tools to prove the second statement of Theorem~\ref{thm:main}.

\noindent\textsl{Proof of Theorem~\ref{thm:main}. }
Let $S$ be an independent set of size $(n-1)!$ in $P(n)$. Assume that
the identity permutation $\sigma$ is in $S$ and let $v_S$ be the
characteristic vector of $S$.

By Lemma~\ref{lem:module}, $v_S$ is in the $\textrm{span}\{\one, h_{i,j}:
i,j \leq n-1 \}$.  We consider two cases, first when $v_S$ is in
$\textrm{span}\{h_{i,j}: i,j \leq n-1\}$ and second when it is not.
 
{\bf case 1.} Assume $v_S \in \textrm{span}\{h_{i,j} :
i,j=1,\dots,n-1\}$, or, equivalently, that $v_S=Hy$ for some vector
$y$.

Since $S$ is an independent set no neighbours of $\sigma$ can be in $S$ and $Ny=0$.
By Lemma~\ref{lem:kernel}, $v_S=Wx$ for some vector $x$.  

For any $i \in \{1, \dots ,n-1\}$ assume the $i$-th entry of the vector $x$ is non-zero.
As $n\geq 3$, there is a permutation $\pi$ with $\pi(i)=i$ and no other fixed
points. This means that the entry in the row corresponding to $\pi$
of $v_S$ must be equal to the $i$-th entry of $x$. Since $v_S$ is a 01-vector,
$x$ must also be a 01-vector.

Further, since $n \geq 4$, for every pair of distinct $i,j \in \{1,
\dots ,n-1\}$ there is a permutation $\pi$ that fixes $i$ and $j$ but
no other points. If the $i$-th and $j$-th entries of $x$ are both
non-zero then the entry in the row corresponding to $\pi$ of $v_S$ is
2. Since $v_S$ must be a 01-vector, there is only one non-zero entry
in $x$. Thus $v_S$ is one of the columns of $W$ and $S = S_{i,i}$ for
some $i \in \{1,\dots,n-1\}$.

{\bf case 2.}  Assume $v_S$ is not in the column space of $H$.
Equivalently, there is some vector $y$ such that $v_S = [H|1]y = Hy_H
+ c\one$ where $y_H$ denotes the vector formed from the first
$(n-1)^2$ entries of $y$ and $c$ is a non-zero constant.

As in case 1, no neighbours of $\sigma$ are in $S$ so $[N|1]y=0$. By
Lemma~\ref{lem:kernelwith1} 
there is a non-zero $c$ such that
\[
y_{[M|1]} = -\frac{c}{(n-2)}(1,1,\dots,1,-(n-2)).
\]
This determines all entries, upto multiplication by a constant, of $y$ except the first $n-1$.

For each $i \leq n-1$ there is a permutation $\pi$ with $\pi(i)=i$ and
no other fixed points.  If $y_i$ is the $i$-th entry of $y$ then
the entry in $v_S$ corresponding to $\pi$ is
\[
y_i + (n-3)\left(-\frac{c}{n-2}\right) +c
\]
which must be either 0 or 1. This implies that
\[
y_i = -\frac{c}{n-2} \textrm{ or } y_i =1 - \frac{c}{n-2}.
\]

Since $n\geq 4$ for any distinct pair $i,j \in \{1,\dots,n-1\}$ there
is a permutation that fixes both $i$ and $j$ and no other points.  If
both $y_i$ and $y_j$ are equal to $1-\frac{c}{n-2}$, then the entry in
the vector $v_S$ which corresponds to this permutation is
\[
2\left(1-\frac{c}{n-2}\right) + (n-4)\left(-\frac{c}{n-2}\right) + c =2
\]
which is a contradiction since $v_S$ is a 01-vector. Thus at most one
of the first $n-1$ entries of $y$ is $1-\frac{c}{n-2}$.

Next, assume that exactly one of the first $n-1$ entries is $1-\frac{c}{n-2}$.
Since $\sigma \in S$ the sum of the first $n-1$ entries of $y$ is 1. But this means that
\[
1-\frac{c}{n-2} + (n-2)\left(-\frac{c}{n-2}\right) + c =1,
\]
which implies that $c=0$, a contradiction.
Hence, all the entries of $y$, except that last, are $-\frac{c}{n-2}$.

Using the fact that the sum of the first $n-1$ entries of $y$ is 1
\[
(n-1)\frac{c}{-(n-2)} + c = 1
\]
which implies that $c=-(n-2)$.

For case 2 there is only one possibility for $y$, this is 
\[
y =(1,1,\dots,-(n-2)). 
\]

Every row in $[H|1]$ that corresponds to a permutation that maps $n$
to $n$ has exactly $(n-1)$ entries equal to one and all other entries
equal to zero. All the other rows has exactly $(n-2)$ entries equal to
one and all other entries equal to zero.  {}From this it follows that $[H|1]y =v_S$ is the characteristic vector of the set $S_{n,n}$.  \qed

\section{Further Work}

We have only considered the simplest version of the Erd\H{o}s-Ko-Rado theorem.
The full version of the Erd\H{o}s-Ko-Rado theorem is concerned with $t$-intersecting subsets.
For an integer $t$, subsets $A,B \subseteq \{1,\dots,n\}$ are \textsl{$t$-intersecting}
if $|A \cap B| \geq t$. 

\begin{theorem}(Erd\H{o}s-Ko-Rado~\cite{MR25:3839})\label{thm:fullekr}
Let $t\leq k \leq n$ be positive integers.
Let $\cA$ be a family of pairwise $t$-intersecting $k$-subsets of $\{1,\dots,n\}$.
There exist a function $f(k,t)$ such that for $n \geq f(n,k)$
\[|\cA| \leq {n-t \choose k-t}.\]
Moreover, a $t$-intersecting family $\cA$ meets this bound if and only if $\cA$ is the collection of all $k$-subsets that contain a fixed $t$-subset.
\end{theorem}

Permutations $\pi, \sigma \in S(n)$ are \textsl{$t$-intersecting} if
\[
|\{i \in \{1,\dots,n\} : \pi(i)=\sigma(i)\}| \geq t.
\]
Again, there is  an obvious family of candidates for the maximum system of $t$-intersecting permutations.
Assume
\[
A = \{(x_i,y_i) : i=1,\dots,t \;\mathrm{ and }\; x_i,y_i \in \{1,\dots,n\} \}
\] 
with $x_i \neq x_j$ and $y_i \neq y_j$ for all $i\neq j$.
Then the family
\[
S_{A} = \{\pi : \pi(x_i) = y_i \;\mathrm{ for \;all }\; (x_i,y_i) \in A\}.
\]
is $t$-intersecting and $|S_A|= (n-t)!$.

Deza and Frankl~\cite{MR0439648} conjecture that
Theorem~\ref{thm:fullekr} can also be extended to families of
$t$-intersecting permutations.

\begin{conj}(Deza and Frankl~\cite{MR0439648}) 
  For $n$ sufficiently
  large, the size of the maximum set of permutations of an $n$-set
  that are pairwise $t$-intersecting is $(n-t)!$.
\end{conj}

Cameron and Ku note that their method cannot be extended to
$t$-intersecting permutations. It is possible that the proof presented in
this paper may be extended as follows.

Define a graph $P_t(n)$ whose vertices are the permutations of an
$n$-set, where two vertices are adjacent if they agree on no more than
$t$ points. Note that $P(n) = P_0(n)$. 

The graph $P_t(n)$ is a sum of all $A_C$ where $C$ is a conjugacy
classes in which the elements have no more than $t$ fixed points.

The graph $P_t(n)$ is vertex transitive so we have that
\[
\alpha(P_t(n))\omega(P_t(n)) \leq n!
\]
This is Equation 3 in Deza and Frankl~\cite{MR0439648}. They also note that if there
exists a sharply $2$-transitive set of permutations of $\{1,\dots,n\}$
(say $PGL(2,n)$) then there is a clique of size $n(n-1)$ and we have
the bound on the $2$-intersecting permutations. 

We conjecture that the shifted characteristic vector of a $2$-intersecting permutation family lies in a union of modules. Define the \textsl{depth} of a partition $\lambda \vdash n$ with $\la=(\la_1,\la_2,\dots,\la_k)$ to be $n-\la_1$.

\begin{conj}
Let $v_S$ be the characteristic vector of a maximum independent set in
$P_1(n)$.  Then the vector $v_S -\frac{|S|}{n!}\one$ lies in the
sum of the modules whose partitions have depth no more than 2. That is
the sum of the following modules
\[
[n],\quad [n-1,1],\quad [n-2,2], \quad [n-2,1,1].
\]
\end{conj}
The dimensions of the sum of these modules and the dimension of the span of $v_{A} -\frac{|S|}{n!}\one$ agree for $n=4,5,6$ where $A=\{(i,j),(k,\ell)\}$.

This conjecture can be generalized to $t$-intersecting permutation systems.
\begin{conj}
Let $v_S$ be the characteristic vector of a maximum independent set in
$P_t(n)$.  Then the vector $v_S -\frac{|S|}{n!}\one$ lies in the
sum of the modules whose partitions have depth no more than $t$.
\end{conj}

Finally, the proof of the Erd\H{o}s-Ko-Rado theorem for permutations given in
this paper is an application of a method that has been used to prove
the Erd\H{o}s-Ko-Rado theorem for set systems and its analogue for
intersecting vector spaces over a finite field. Another direction for
this work is to apply this method to other objects such as perfect
matchings and uniform partitions with a plan of developing a more
general theory of Erd\H{o}s-Ko-Rado theorems.

\end{document}